\documentclass[preprint]{elsart}
\usepackage{amsfonts}
\usepackage{amsmath}
\usepackage{graphicx}
\usepackage{tikz}
\usetikzlibrary{trees}
\usetikzlibrary{arrows,shapes,snakes,automata,backgrounds,petri}

\usepackage{epstopdf}
\begin{document}
\begin{frontmatter}

\title{Graph-theoretic criteria for injectivity and unique equilibria in general chemical reaction systems}

\author{Murad Banaji\thanksref{cor1}\thanksref{ref1} }
\and
\author{Gheorghe Craciun\thanksref{ref2}}

\address[ref1]{Department of Mathematics, University College London, Gower Street, London WC1E 6BT, UK and Department of Biological Sciences, University of Essex, Wivenhoe Park, Colchester, CO4 3SQ, UK. }
\address[ref2]{Department of Mathematics and Department of Biomolecular Chemistry, University of Wisconsin, Madison, WI 53706-1388.}
\thanks[cor1]{Corresponding author. m.banaji@ucl.ac.uk.}

\begin{abstract}
In this paper we discuss the question of how to decide when a general chemical reaction system is incapable of admitting multiple equilibria, regardless of parameter values such as reaction rate constants, and regardless of the type of chemical kinetics, such as mass-action kinetics, Michaelis-Menten kinetics, etc. Our results relate previously described linear algebraic and graph-theoretic conditions for injectivity of chemical reaction systems. After developing a translation between the two formalisms, we show that a graph-theoretic test developed earlier in the context of systems with mass action kinetics, can be applied to reaction systems with arbitrary kinetics. The test, which is easy to implement algorithmically, and can often be decided without the need for any computation, rules out the possibility of multiple equilibria for the systems in question. 
\end{abstract}


\end{frontmatter}

\section{Introductory material}

There is increasing interest in methodologies for drawing conclusions about the {\em dynamics} of a chemical reaction network based only on the {\em network structure}, i.e., with limited or absent knowledge of the kinetics. Early work in this direction is exemplified by \cite{hornjackson,feinberg0,horn,feinberg}, with more recent strands including discussions of monotonicity \cite{kunzesiegelmathchem,leenheer,banajidynsys}, and discussions of injectivity \cite{craciun,craciun1,banajiSIAM}. Although the reaction systems discussed in examples are often assumed to have mass action kinetics, an important feature of \cite{leenheer,banajidynsys,banajiSIAM} was that they described criteria which applied for more general chemical kinetics. The aim of this paper is to extend graph-theoretic results, which were developed in \cite{craciun1} in the context of mass action kinetics, to the case of general chemical kinetics. These graph-theoretic criteria are more restrictive than the matrix-theoretic results in \cite{banajiSIAM}, but are more intuitive, and give rise to conditions which are less expensive computationally, and are often easy to check by hand. 

{\bf Dynamical systems derived from chemical reaction networks.} A chemical reaction system in which $n$ reactants participate in $m$ reactions has dynamics governed by the ordinary differential equation
\begin{equation}
\dot x = Sv(x), \label{basic}
\end{equation}
where $x = [x_1, \ldots, x_n]^T$ is the nonnegative $n$-vector of reactant concentrations, $v = [v_1, \ldots, v_m]^T$ is the $m$-vector of reaction rates and $S$ is the $n \times m$ stoichiometric matrix. Arbitrary orderings can be chosen on the sets of substrates and reactions. Further, $S$ is only defined up to an arbitrary re-signing of its columns, equivalent to a switching of the left and right hand sides of a reaction. It is trivial that all results here are independent of the orders chosen on substrates and reactions. We will also confirm below that they are independent of the signing of columns of $S$. System~(\ref{basic}) defines a dynamical system on the nonnegative orthant in $\mathbb{R}^n$. With the additional assumption that all substances may have some inflow (which is allowed to be zero) and some outflow which increases with concentration, we obtain the related system
\begin{equation}
\dot x = K + Sv(x) - Q(x). \label{basicout}
\end{equation}
Here $K$ is a nonnegative vector representing the inflows and the diagonal function $Q(x) = [Q_1(x_1), \ldots, Q_n(x_n)]^T$ represents the outflows, and is assumed to satisfy $\frac{\partial Q_i}{\partial x_i} > 0$ for each $i$. The system has Jacobian $J = SV(x) - D(x)$ where the $m \times n$ matrix $V(x)$ is defined by $V_{ij}(x) \equiv \frac{\partial v_i}{\partial x_j}$, and the diagonal matrix $D(x)$ is defined by $D_{ii}(x) \equiv \frac{\partial Q_i}{\partial x_i}$. It is notationally convenient to omit the explicit dependence on $x$, and write $J=SV-D$.

{\bf Previous results on injectivity and the main result of this paper.} In \cite{banajiSIAM} and \cite{banajidynsys} a reaction system was termed ``nonautocatalytic'' (NAC for short) if $S$ and $V^T$ have opposite sign structures in the following sense: $S_{ij}V_{ji} \leq 0$ for all $i$ and $j$, and $S_{ij} = 0 \Rightarrow V_{ji} = 0$. These conditions are naturally fulfilled provided that no substrate occurs on both sides of a reaction (either with the same or with a changed stoichiometry). Here we refer to systems where the above two conditions are fulfilled, as {\bf N1C} reaction systems, in order to emphasise that the conditions only rule out {\em one-step} catalysis. Note that the  {\bf N1C} condition is not very restrictive for realistic biochemical reaction networks: for example, all the enzymatic reactions considered in \cite{craciun3} satisfy this condition.

A square matrix $M$ is {\bf sign nonsingular} if all matrices with the same sign pattern as $M$ are nonsingular \cite{brualdi}. In \cite{banajiSIAM} it was shown that System~(\ref{basicout}) is injective, and hence incapable of multiple equilibria, provided that the stoichiometric matrix $S$ is {\bf strongly sign determined} ({\bf SSD}), i.e., all square submatrices of $S$ are either sign nonsingular or singular. This led to a characterization of injectivity based entirely on a computation on the stoichiometric matrix. On the other hand, following theory developed in \cite{craciun}, a signed, labelled, bipartite multigraph termed the Species-Reaction graph or {\bf SR graph} was constructed in \cite{craciun1} and used to make claims about the nonexistence of multiple equilibria in systems of chemical reactions with mass-action kinetics. This time, rather than a matrix computation, a certain graph-theoretic condition, sometimes checkable by observation alone without the need for any computation, was shown to be sufficient to guarantee the absence of multiple equilibria. Here we will combine these ideas, and in fact show that the graph-theoretic condition in \cite{craciun1} suffices to guarantee the absence of multiple equilibria for {\bf N1C} systems with arbitrary kinetics. 

We will define an SR graph for an arbitrary {\bf N1C} reaction system solely via the stoichiometric matrix $S$ of the system. We will then present a condition on this graph which will be termed Condition~($*$). Our main result is the following:
\begin{thm}
\label{suff}
Consider the SR graph $G$ of an {\bf N1C} reaction system with stoichiometric matrix $S$. Assume that Condition~($*$) is fulfilled. Then $S$ is an SSD matrix. 
\end{thm}

As a consequence, if Condition~($*$) holds for the SR graph associated with an N1C reaction system, then System~(\ref{basicout}), i.e. the system with outflows, is injective, and hence incapable of multiple equilibria. Along the way to proving Theorem~\ref{suff} we will prove several auxiliary results relating computations on matrices to computations on graphs. Of these, Corollary~\ref{SNS2} is of independent interest. 

\section{The SR graph}

The SR graph, introduced in \cite{craciun1}, is a bipartite graph with $n$ substrate vertices or {\bf S-vertices}, $m$ reaction vertices or {\bf R-vertices}, and an edge between S-vertex $i$ and R-vertex $j$ iff substrate $i$ participates in reaction $j$. For an {\bf N1C} reaction system there is a one-to-one correspondence between edges in the SR graph and nonzero entries in the stoichiometric matrix $S$: clearly $S_{ij} \not = 0$ implies that substrate $i$ participates in reaction $j$, and given the {\bf N1C} condition it can participate on one side of the reaction only, giving one edge between S-vertex $i$ and R-vertex $j$ (in the general case multiple edges are allowed). On the other hand, given the {\bf N1C} condition $S_{ij} = 0$ implies that substrate $i$ does not participate in reaction $j$, and hence that there is no edge connecting S-vertex $i$ and R-vertex $j$.

In this paper, for convenience, we make slight modifications to the labelling of the SR graph in \cite{craciun1}. An edge $e_{ij}$ (where $ i \in \{1, \ldots, n\}$ and $j \in \{1, \ldots, m\}$) will be taken to mean an edge between the $i$th S-vertex and the $j$th R-vertex. Such an edge exists if and only if $S_{ij} \not = 0$. Rather than labelling edge $e_{ij}$ with the complex label associated with substrate $i$ and reaction $j$ as done in \cite{craciun1}, we simply give it a sign, so that  $\mathrm{sign}(e_{ij}) = \mathrm{sign}(S_{ij})$. Since $S$ is only defined upto an arbitrary signing its columns, the signing of edges in the SR graph is not unique. However, as we will show in Lemma~\ref{resigning} after developing some further ideas, all results here are independent of the choice of signing. Whether labelled with the complex labels or signed, the important thing is that the SR graph contains information on whether substrates occur on the same side of a reaction or on opposite sides of a reaction.

Given the one-to-one correspondence
\[
S_{ij} \Leftrightarrow e_{ij}
\]
when $S_{ij} \not = 0$ it is convenient to allow a slight abuse of notation and refer to $S_{ij}$ as an edge in $G$. Two edges $S_{i_1j_1}$ and $S_{i_2j_2}$ share an S-vertex when $i_1 = i_2$, and they share an R-vertex when $j_1 = j_2$. 

{\bf Remark.} As a consequence of our formal redefinition, rather than being associated directly with chemical reaction networks, SR graphs are now associated directly with real matrices. 

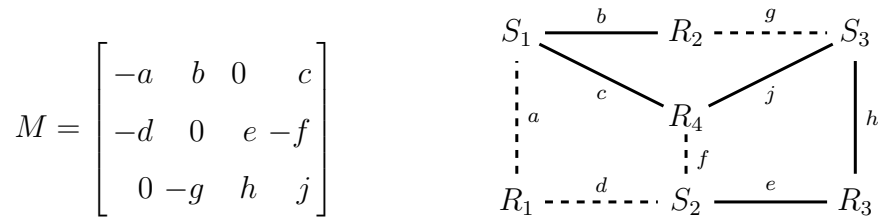
\begin{figure}[h]
\begin{minipage}{0.4\textwidth}
\[
M = \left[
\begin{array}{rrrr}
-a & b & \,\,\,0\,\, & c\\-d &  0 & \,\,e & -f\\ 0 &  -g & \,\,h & j
\end{array}
\right]
\]

\end{minipage}
\hfill
\begin{minipage}{0.5\textwidth}
\begin{tikzpicture}[domain=0:4,scale=0.75]

\node at (1,4) {$S_1$};
\node at (1,1) {$R_1$};
\node at (4,4) {$R_2$};
\node at (4,1) {$S_2$};
\node at (7,4) {$S_3$};
\node at (7,1) {$R_3$};

\node at (4,2.5) {$R_4$};

\draw[-, dashed, line width=0.04cm] (1,1.5) -- (1, 3.5);
\draw[-, line width=0.04cm] (7,1.5) -- (7, 3.5);

\draw[-, dashed, line width=0.04cm] (1.5,1) -- (3.5, 1);
\draw[-, line width=0.04cm] (4.5,1) -- (6.5, 1);
\draw[-, line width=0.04cm] (1.5,4) -- (3.5, 4);
\draw[-, dashed, line width=0.04cm] (4.5,4) -- (6.5, 4);

\draw[-, dashed, line width=0.04cm] (4,1.5) -- (4, 2.2);

\draw[-, line width=0.04cm] (3.6,2.7) -- (1.4, 3.8);
\draw[-, line width=0.04cm] (4.4,2.7) -- (6.6, 3.8);

\node at (2.5, 4.3) {$\scriptstyle{b}$};
\node at (2.5, 1.3) {$\scriptstyle{d}$};

\node at (5.5, 4.3) {$\scriptstyle{g}$};
\node at (5.5, 1.3) {$\scriptstyle{e}$};

\node at (1.3, 2.5) {$\scriptstyle{a}$};
\node at (7.3, 2.5) {$\scriptstyle{h}$};

\node at (4.3, 1.75) {$\scriptstyle{f}$};
\node at (2.5, 2.9) {$\scriptstyle{c}$};
\node at (5.5, 2.9) {$\scriptstyle{j}$};

\end{tikzpicture}
\end{minipage}
\caption{\label{exSR} {\em Left.} A matrix $M$. Assume that $a, b, c,d,e,f,g,h,j > 0$. {\em Right.} The corresponding SR graph. S-vertices have been labelled $S_1, S_2, S_3$ while R-vertices have been labelled $R_1, R_2, R_3, R_4$. Positive edges are bold lines while negative edges are dashed lines. These conventions will be followed in drawings of SR graphs.}
\end{figure}

{\bf Example.} An example of a matrix and its corresponding SR graph is shown in Figure~\ref{exSR}.

{\bf Terminology and definitions.} Given a graph $G$, we will say that some vertex/edge lies in $G$ meaning that it lies in the vertex/edge-set of $G$. We will generally define subgraphs by their edge-sets, i.e., writing $E = \{e_1, \ldots, e_k\}$ where $e_i$ are edges in $G$, will mean that $E$ is the subgraph of $G$ containing exactly edges $e_1, \ldots, e_k$ and the vertices on which they are incident. The terms {\bf path} and {\bf cycle} will be taken to mean open/closed simple walks in $G$ respectively and are particular examples of subgraphs of $G$. In the usual way, the size of a subgraph $E$ is the number of edges in $E$, written $|E|$. When $E$ is a cycle or a path this will also be called the length of $E$. Because of the bipartite nature of SR graphs, all cycles are of even length including alternate S- and R-vertices. If two paths or cycles are edge-disjoint they may still share some vertices. When they share no vertices they will be termed vertex-disjoint. A path between an S-vertex and an R-vertex is called an {\bf S-to-R path}. Note that S-to-R paths are of odd length.

Cycles in an SR graph have a natural {\bf parity} -- they are either odd or even. We use a definition equivalent to that in \cite{craciun1}, but using the labelling defined above.

Consider a matrix $S$ and the corresponding SR graph $G$. Let $E$ be any subgraph of $G$. The sign of an edge $e$ in $G$ has already been defined above, from which we can define the {\bf sign} of $E$ to be 
\[
\mathrm{sign}(E) = \prod_{e \in E}\mathrm{sign}\left(e\right).
\]
When $|E|$ is even, we define the parity of $E$ to be
\[
P(E) = (-1)^{|E|/2}\mathrm{sign}(E).
\]
Since cycles are always of even length, the parity of a cycle is always defined. A cycle $C$ will be termed an {\bf e-cycle} if $P(C) = 1$ and an {\bf o-cycle} if $P(C) = -1$. Note that by these definitions, for an e-cycle $C$ we have
\[
\mathrm{sign}(C) = (-1)^{|C|/2},
\] 
and similarly for an o-cycle $C$ we have
\[
\mathrm{sign}(C) = (-1)^{|C|/2 - 1}.
\] 

We define the {\bf value} of edge $e_{ij}$ to be $\mathrm{val}(e_{ij}) = |S_{ij}|$, and for a subgraph $E$, 
\[
\mathrm{val}(E) = \prod_{e \in E}\mathrm{val}(e)
\]
When $C$ is a cycle containing edges $e_1, e_2, \ldots, e_{2r}$ such that $e_i$ and $e_{i+1 \mod 2r}$ are adjacent for each $i =1,\ldots,2r$, we can define
\[
\mathrm{stoich}(C) = \left|\prod_{i = 1}^{r}\mathrm{val}(e_{2i-1}) - \prod_{i = 1}^{r}\mathrm{val}(e_{2i})\right|.
\]
Note that this definition is independent of the starting point chosen on the cycle. A cycle with $\mathrm{stoich}(C) = 0$ is termed an {\bf s-cycle}. This definition is equivalent to the definition in \cite{craciun1}. 

The intersection of two cycles in an SR graph can be divided into a set of vertex-disjoint paths. We say that two cycles have {\bf S-to-R intersection}, if each component of their intersection is an S-to-R path, i.e. a path between an S-vertex and an R-vertex.

We finish this section with a lemma confirming that re-signing the columns of the stoichiometric matrix does not alter the nature of cycles in an SR graph. 
\begin{lem}
\label{resigning}
Consider a matrix $S$ and the associated SR graph $G$. Now consider any re-signed version of $S$, say $S^{'}$ with corresponding SR graph $G^{'}$. Cycles in $G$ are in one-to-one correspondence with those in $G^{'}$. e-cycles (o-cycles) in $G$ correspond to e-cycles (o-cycles) in $G^{'}$. s-cycles in $G$ correspond to s-cycles in $G^{'}$.
\end{lem}
\begin{pf}
Ignoring the signs on edges, clearly $G$ and $G^{'}$ are identical graphs, so cycles in $G$ and $G^{'}$ are in one-to-one correspondence. Consider some cycle $C \in G$ and the corresponding cycle $C^{'} \in G^{'}$. Re-signing column $j$ of $S$ means re-signing all edges incident on R-vertex $j$ in $G^{'}$. But clearly $C$, and hence $C^{'}$, contains either 2 or 0 edges incident on R-vertex $j$. So the re-signing causes no change in $\mathrm{sign}(C)$, i.e., $\mathrm{sign}(C) = \mathrm{sign}(C^{'})$, so that e-cycles and o-cycles are preserved. Finally, re-signing a column does not alter any of the values of edges, leaving $\mathrm{stoich}(C)$ unchanged. Thus the re-signing does not affect the property of a cycle being an s-cycle. \qquad \qed
\end{pf}

As a final note, it is obvious that S-to-R paths remain so after a re-signing of the graph.

\section{Determinants, permutations and cycles}
\label{detperm}

Consider some $n \times m$ matrix $S$ and the associated SR graph $G$. 

{\bf Notation.} $S(\gamma|\delta)$ will refer to the submatrix of $S$ with rows from some set $\gamma\subset\{1, \ldots, n\}$ and columns from some set $\delta\subset\{1, \ldots, m\}$. If $S(\gamma|\delta)$ is square, then its determinant will be written $S[\gamma|\delta]$. Each submatrix $S(\gamma|\delta)$ corresponds to a subgraph of $G$ which we will term $G(\gamma|\delta)$. 

{\bf Determinant expansions and term subgraphs.} 
Consider any sets $\gamma = \{\gamma_1, \gamma_2, \ldots, \gamma_k\} \subset \{1, \ldots, n\}$ and $\delta = \{\delta_1, \delta_2, \ldots, \delta_k\} \subset \{1, \ldots, m\}$, choosing $\gamma_1 < \gamma_2 < \ldots < \gamma_k$ and $\delta_1 < \delta_2 < \ldots < \delta_k$, so that both sets have a natural ordering. Consider the square submatrix $S(\gamma|\delta)$ of $S$. Any permutation $\alpha$ of $\delta$ gives a term in the expansion of $S[\gamma|\delta]$ of the form $T_\alpha = P(\alpha)S_{\gamma_1 \alpha_1}\cdots S_{\gamma_k \alpha_k}$, where $P(\alpha) = -1$ if $\alpha$ is an odd permutation, and $P(\alpha) = 1$ otherwise. Note that $T_\alpha$ is simply a real number. Thus terms in the determinant expansion of $S(\gamma|\delta)$ are in one-to-one correspondence with permutations of $\delta$. 

If $T_\alpha \not = 0$, then $\alpha$ can also be identified with a subgraph of $G(\gamma|\delta)$, 
\[
E_\alpha = \{S_{\gamma_1\alpha_1}, \ldots, S_{\gamma_k\alpha_k} \}.
\]
Each S-vertex chosen from $\gamma$ and each R-vertex from $\delta$ occurs in exactly one edge in $E_\alpha$. This follows since each member of $\{\gamma_1, \ldots, \gamma_k\}$ occurs exactly once as a first subscript in $S_{\gamma_1 \alpha_1},\ldots, S_{\gamma_k \alpha_k}$ and similarly each of $\{\delta_1, \ldots, \delta_k\}$ occurs exactly once as a second subscript in this expression. As a result no two edges in $E_\alpha$ share a vertex, and $E_\alpha$ can contain no cycles. We will refer to a subgraph of $G(\gamma|\delta)$ with these properties as a {\bf term subgraph} of $G(\gamma|\delta)$. Clearly term subgraphs are in one-to-one correspondence with nonzero terms in the expansion of $S[\gamma|\delta]$. 

Permutations of a fixed set $\delta$ form a group, and so it makes sense to talk about the operations of composition and inversion. Given permutations $\alpha, \beta$ we will often be interested in the permutation $\beta\circ \alpha^{-1}$, which takes $\alpha$ to $\beta$. Permutations can be written as products of {\em disjoint} cycles. A {\bf nontrivial cycle} will refer to a cycle of length greater than 1. Below we will show that there is a close relationship between cycles in the decomposition of a permutation, and cycles in the SR graph. Throughout this paper, the word ``cycle'' has two distinct meanings: either a special kind of subgraph in an SR graph, or a special kind of permutation. In general the meaning will be clear from the context. 

Several key constructions in this paper rely on taking two term subgraphs corresponding to two distinct permutations of $\delta$, say $\alpha$ and $\beta$, and looking at the structure of their union $E_\alpha \cup E_\beta$. Consider this union: the S-vertex $\gamma_i$ occurs in exactly one edge in $E_\alpha \cup E_\beta$ if $\alpha_i = \beta_i$, and exactly two edges in $E_\alpha \cup E_\beta$ otherwise. On the other hand, assume $\delta_j = \alpha_i = \beta_k$: if $i=k$, then the R-vertex $\delta_j$ occurs in exactly one edge in $E_\alpha \cup E_\beta$ while if $i \not = k$ then it occurs in exactly two edges in $E_\alpha \cup E_\beta$. Thus $E_\alpha \cup E_\beta$ consists of a set of components, each of which is either i) an isolated edge corresponding to a trivial cycle in $\beta\circ \alpha^{-1}$, or ii) a cycle corresponding to a nontrivial cycle in $\beta\circ \alpha^{-1}$. The explicit construction is carried out in Lemma~\ref{cycles} below. Any given vertex in $E_\alpha \cup E_\beta$ has exactly one edge from $E_\alpha \cup E_\beta$ incident on it if it corresponds to a trivial cycle in $\beta\circ\alpha^{-1}$, and exactly two edges incident on it if it corresponds to a nontrivial cycle. In a similar way, given a set of $k$ permutations $\alpha_1, \ldots, \alpha_k$, each vertex in $\cup_{i=1}^k E_{\alpha_i}$ has between $1$ and $k$ edges from $\cup_{i=1}^k E_{\alpha_i}$ incident on it. 

{\bf Notation.} Given a set $\gamma = \{\gamma_1, \gamma_2, \ldots, \gamma_k\}$, we write $\gamma = [\gamma_1, \gamma_2, \ldots, \gamma_k]$, when it is important to stress that $\gamma$ is an ordered set. When a permutation is written as a product of cycles, we use round brackets to denote these cycles. These may include or exclude trivial cycles. 

{\bf Example.} Let 
\[
\delta = [\delta_1, \delta_2, \delta_3, \delta_4], \qquad \alpha = [\delta_1, \delta_3, \delta_4, \delta_2], \qquad \beta = [\delta_2, \delta_4, \delta_3, \delta_1].
\]
Written as products of disjoint cycles, $\alpha = (\delta_1)(\delta_2, \delta_3, \delta_4)$, $\beta = (\delta_3)(\delta_1, \delta_2, \delta_4)$ and $\beta\circ\alpha^{-1} = (\delta_1 , \delta_2)(\delta_3, \delta_4)$. Finally $P(\alpha) = P(\beta) = +1$ and $P(\beta\circ\alpha^{-1}) = P(\beta)P(\alpha) = +1$. 

{\bf Notation.} Given that cycles are the fundamental objects in the theory being developed in this paper, when we have an index $k$ which is known to belong to a set $\{1, \ldots, r\}$, counting is always done on a circle of size $r$, so that $k+p$ means $(k+p-1 \mod r) + 1$. Adopting this convention avoids lengthy subscripts.

The next lemma illustrates the relationship between cycles in permutations and cycles in SR graphs, and also begins the process of linking statements about SR graphs and stoichiometric matrices. 
\begin{lem}
\label{cycles}
Consider a matrix $S$ and the associated SR graph $G$. Assume that there is some square submatrix $S(\gamma|\delta)$ and permutations $\alpha$ and $\beta$ of $\delta$ such that $T_\alpha$ and $T_\beta$ are nonzero terms in $S[\gamma|\delta]$. Then corresponding to each nontrivial cycle in the decomposition $\beta\circ\alpha^{-1}$ there is a cycle in $G$. In particular $G$ contains at least one cycle. 
\end{lem}

\begin{pf}
The two nonzero terms can be written explicitly as:
\[
T_\alpha = P(\alpha)S_{\gamma_1 \alpha_1}\cdots S_{\gamma_k \alpha_k}, \qquad T_\beta = P(\beta)S_{\gamma_1 \beta_1}\cdots S_{\gamma_k \beta_k}\,.
\]
Write  $\beta \circ \alpha^{-1}$ as the product of disjoint cycles. Since $\beta \not = \alpha$ this product contains at least one nontrivial cycle. Consider such a cycle of length $r$ ($2 \leq r \leq k$), $(\delta_{b(1)}, \ldots, \delta_{b(r)})$. For each $j \in \{1, \ldots, r\}$ define the index $a(j)$ by $\alpha_{a(j)} = \delta_{b(j)}$. (Note that since $\alpha$ is a permutation, $j_1 \not = j_2 \Leftrightarrow a(j_1) \not = a(j_2)$.) The existence of the cycle means that for each $j\in \{1, \ldots, r\}$, $S_{\gamma_{a(j)}\delta_{b(j)}} \in E_\alpha$ and $S_{\gamma_{a(j)}\delta_{b(j+1)}} \in E_\beta$. Writing these as \[\{S_{\gamma_{a(1)}\delta_{b(1)}}, S_{\gamma_{a(1)}\delta_{b(2)}}, S_{\gamma_{a(2)}\delta_{b(2)}}, \ldots, S_{\gamma_{a(r)}\delta_{b(r)}}, S_{\gamma_{a(r)}\delta_{b(1)}}\}\] makes it clear that they form a cycle of length $2r$ in $G$. 
\qquad \qquad \qed
\end{pf}

The construction in Lemma~\ref{cycles} will be used frequently -- i.e., given two different permutations $\alpha$ and $\beta$ of a set $\delta$, cycles in $\beta\circ\alpha^{-1}$ will be used to infer the existence of index sets $\{a(j)\}$ and $\{b(j)\}$ and corresponding cycles in the SR graph. When using the construction, for notational brevity we will write $a_1$ for $\gamma_{a(1)}$ and $b_1$ for $\delta_{b(1)}$.

{\bf Example.} Let 
\[
\begin{array}{cclccr}
\gamma & = & [\gamma_1, & \gamma_2, & \gamma_3, & \gamma_4]\\
\delta & = & [\delta_1, & \delta_2, & \delta_3, & \delta_4]\\
\alpha & = & [\delta_1, & \delta_3, & \delta_4, & \delta_2]\\
\beta & = & [\delta_2, & \delta_4, & \delta_3, & \delta_1]
\end{array}
\]
Let $S(\gamma|\delta)$ be a submatrix of a matrix $S$ and $G(\gamma|\delta)$ the associated subgraph. $\alpha$ and $\beta$ define the term subgraphs of $G(\gamma|\delta)$: $E_\alpha = \{S_{\gamma_1\delta_1}, S_{\gamma_2\delta_3}, S_{\gamma_3\delta_4}, S_{\gamma_4\delta_2}\}$ and $E_\beta = \{S_{\gamma_1\delta_2}, S_{\gamma_2\delta_4}, S_{\gamma_3\delta_3}, S_{\gamma_4\delta_1}\}$. From the previous example, $\beta\circ\alpha^{-1} = (\delta_1 , \delta_2)(\delta_3, \delta_4)$, so there are two cycles $C_1$ and $C_2$ in $E_\alpha \cup E_\beta$. Written out explicitly, these cycles are \[
C_1 = \{S_{\gamma_1\delta_1}, S_{\gamma_1\delta_2}, S_{\gamma_4\delta_2}, S_{\gamma_4\delta_1}\}\quad \mbox{and} \quad C_2 = \{S_{\gamma_2\delta_3}, S_{\gamma_2\delta_4}, S_{\gamma_3\delta_4}, S_{\gamma_3\delta_3}\}.
\]

So far we have focussed on constructing cycles in an SR graph from pairs of nonzero terms in a determinant. However the reverse is also important -- inferring the existence of pairs of nonzero terms in a determinant from structures in a graph. The basic operation which allows us to do this is a particular splitting of a cycle. Any cycle
\[
C = \{S_{a_1 b_1}, S_{a_1 b_2}, S_{a_2 b_2}, \ldots, S_{a_r b_r}, S_{a_r b_1}\}
\]
can be uniquely partitioned into two vertex-disjoint subgraphs of equal size:
\begin{equation}
\label{cycdecomp}
C = \{S_{a_1 b_1}, S_{a_2 b_2}, \ldots, S_{a_r b_r}\} \cup \{S_{a_1 b_2}, S_{a_2 b_3}, \ldots, S_{a_r b_1}\}.
\end{equation}
We will call this a {\bf disconnecting partition} of $C$. With this notion, confirming if a cycle $C$ is an $s$-cycle now involves:
\begin{enumerate}
\item Constructing a disconnecting partition of $C$ into $\{C_1, C_2\}$,
\item Confirming that $\mathrm{val(C_1)} = \mathrm{val(C_2)}$.
\end{enumerate}

\section{Preliminary lemmas}
With the machinery set up above, we are ready to prove some lemmas. The idea of these lemmas is that each of them is quite brief and reusable, so that the proofs of the main results in this paper become simpler.

The first lemma gives us a basic way of checking whether a permutation, written as a product of cycles, is even or odd. 
\begin{lem}
\label{permsigns}
Consider a permutation $\alpha$ written as a product of cycles from some set $\mathcal{C}$. Let $\theta = \cup_{c \in \mathcal{C}}c$. Then 
\[
P(\alpha) = (-1)^{\left(\left|\theta\right| - \left|\mathcal{C}\right|\right)}
\]
i.e., $\alpha$ is an even or odd permutation according to whether the total number of elements in cycles, minus the number of cycles is even or odd. 
\end{lem}

\begin{pf}
This follows by writing any permutation as the product of disjoint cycles and noting the elementary result that a $k$-cycle is an even permutation if $k$ is odd and vice versa.
\qquad \qed
\end{pf}

In the following lemmas we will pass without comment between talking about terms or sets of terms in submatrices of a matrix, and subgraphs or cycles in the corresponding SR graph. The next formula gives us a way of deciding whether two terms in a determinant expansion have the same sign. 
\begin{lem}
\label{prodforumla}
Consider any square submatrix $S(\gamma|\delta)$ of a matrix $S$. Consider any two nonzero terms $T_\alpha$ and $T_\beta$ in the determinant expansion of $S[\gamma|\delta]$ corresponding to permutations $\alpha$ and $\beta$ of $\delta$. Then
\begin{equation}
\label{prodform}
\mathrm{sign}(T_\alpha T_\beta)  = (-1)^{|\mathcal{C}_e|}
\end{equation}
where $|\mathcal{C}_e|$ is the number of e-cycles in $E_\alpha \cup E_\beta$.
\end{lem}

\begin{pf}
Let $k=|\gamma|=|\delta|$. By definition
\[
T_\alpha T_\beta = P(\alpha)P(\beta)\prod_{i =1}^k S_{\gamma_i \alpha_i}S_{\gamma_i \beta_i}\,.
\]
Let $\theta$ be the set of indices for which $\alpha_i \not = \beta_i$. We can write
\[
T_\alpha T_\beta = P(\alpha)P(\beta)\prod_{i \in \{1, \ldots, k\}\backslash\theta} S_{\gamma_i \alpha_i}S_{\gamma_i \beta_i}\prod_{i \in \theta} S_{\gamma_i \alpha_i}S_{\gamma_i \beta_i}\,.
\]

When $i \in \{1, \ldots, k\}\backslash\theta$, $S_{\gamma_i \alpha_i}S_{\gamma_i \beta_i} = S_{\gamma_i \alpha_i}^2 > 0$. So
\[
\mathrm{sign}(T_\alpha T_\beta) = P(\alpha)P(\beta)\,\,\mathrm{sign}\left(\prod_{i \in \theta} S_{\gamma_i \alpha_i}S_{\gamma_i \beta_i}\right)\,.
\]

Write $\beta \circ \alpha^{-1}$ as a product of disjoint nontrivial cycles, and identify each of these with a cycle in $G(\gamma|\delta)$. Let the set of o-cycles in this set be $\mathcal{C}_o$ and the set of e-cycles be $\mathcal{C}_e$, with  $\mathcal{C} =  \mathcal{C}_o \cup \mathcal{C}_e$. Associate with each cycle $c \in \mathcal{C}_o \cup \mathcal{C}_e$ the corresponding index set $\tilde{c}$, i.e., $i \in \tilde{c} \Leftrightarrow S_{\gamma_i\alpha_i}, S_{\gamma_i\beta_i} \in c$. Thus corresponding to the sets $ \mathcal{C}_o$ and $ \mathcal{C}_e$ are the sets of index sets $\tilde{\mathcal{C}_o}$ and $\tilde{\mathcal{C}_e}$. Since any two cycles are edge-disjoint, $\tilde{\mathcal{C}_o} \cup \tilde{\mathcal{C}_e}$ is a partition of $\theta$, and we can define
\[
\theta_o \equiv \bigcup_{\tilde c \in \tilde{\mathcal{C}_o}}\tilde c,\quad\theta_e \equiv \bigcup_{\tilde c \in \tilde{\mathcal{C}_e}}\tilde c\quad\mbox{with}\quad |\theta_o| = \sum_{\tilde{c} \in \tilde{\mathcal{C}_o}}|\tilde{c}|,\quad |\theta_e| = \sum_{\tilde{c} \in \tilde{\mathcal{C}_e}}|\tilde{c}|\,.
\] 
Clearly $\theta = \theta_o \cup \theta_e$. We can write
\begin{eqnarray}
\prod_{i \in \theta} S_{\gamma_i\alpha_i}S_{\gamma_i\beta_i} & = & \left(\prod_{i \in \theta_o} S_{\gamma_i\alpha_i}S_{\gamma_i\beta_i}\right)\left(\prod_{i \in \theta_e} S_{\gamma_i\alpha_i}S_{\gamma_i\beta_i}\right) \nonumber\\
&  = &\left(\prod_{\tilde{c} \in \tilde{\mathcal{C}_o}}\prod_{i \in \tilde{c}} S_{\gamma_i\alpha_i}S_{\gamma_i\beta_i}\right)\left(\prod_{\tilde{c} \in \tilde{\mathcal{C}_e}}\prod_{i \in \tilde{c}} S_{\gamma_i\alpha_i}S_{\gamma_i\beta_i}\right) \,.\nonumber
\end{eqnarray}

So
\begin{eqnarray}
\hspace{-0.7cm}\mathrm{sign}(T_\alpha T_\beta) & = & P(\alpha)P(\beta)\left(\prod_{\tilde{c} \in \tilde{\mathcal{C}_o}}\mathrm{sign}\left(\prod_{i \in \tilde{c}} S_{\gamma_i\alpha_i}S_{\gamma_i\beta_i}\right)\right)\left(\prod_{\tilde{c} \in \tilde{\mathcal{C}_e}}\mathrm{sign}\left(\prod_{i \in \tilde{c}} S_{\gamma_i\alpha_i}S_{\gamma_i\beta_i}\right)\right) \nonumber \\
& = & P(\alpha)P(\beta)\left(\prod_{c \in \mathcal{C}_o}(-1)^{|c|-1}\right)\left(\prod_{c \in \mathcal{C}_e}(-1)^{|c|}\right) \nonumber\\
& = & P(\alpha)P(\beta)(-1)^{|\theta_o| + |\theta_e|- |\mathcal{C}_o|} \nonumber\\
& = & P(\alpha)P(\beta)(-1)^{|\theta|- |\mathcal{C}_o|} \,. \nonumber
\end{eqnarray}
Applying Lemma~\ref{permsigns} to $\beta\circ\alpha^{-1}$ gives us that 
\[
P(\alpha)P(\beta) = P(\beta\circ\alpha^{-1}) = (-1)^{\left|\theta\right| - \left|\mathcal{C}\right|}\,,
\]
so that
\[
\mathrm{sign}(T_\alpha T_\beta) = (-1)^{\left|\theta\right| - \left|\mathcal{C}\right|}(-1)^{|\theta|- |\mathcal{C}_o|} = (-1)^{2|\theta| - |\mathcal{C}| - |\mathcal{C}_o|} = (-1)^{|\mathcal{C}| + |\mathcal{C}_o|} = (-1)^{|\mathcal{C}_e|}\,.
\]
This proves the result. \qquad \qed
\end{pf}

A corollary of the previous lemma is that when the intersection of two term subgraphs contains only o-cycles then the two corresponding terms have the same sign.
\begin{lem}
\label{allocycles}
Consider any square submatrix $S(\gamma|\delta)$ of a matrix $S$. Let $\alpha$ and $\beta$ be permutations of $\delta$ such that $T_\alpha$ and $T_\beta$ are nonzero terms in the determinant expansion of $S[\gamma|\delta]$. If all cycles in $E_\alpha\cup E_\beta$ are o-cycles, then $\mathrm{sign}(T_\alpha) = \mathrm{sign}(T_\beta)$.  
\end{lem}
\begin{pf}
Since in this case, there are no e-cycles in $E_\alpha \cup E_\beta$, applying Eq.~(\ref{prodform}) immediately gives
\[
\mathrm{sign}(T_\alpha T_\beta)  = (-1)^{|\mathcal{C}_e|} = (-1)^{0} = 1\,.
\]
Thus $\mathrm{sign}(T_\alpha) = \mathrm{sign}(T_\beta)$. 
\qquad \qed
\end{pf}

{\bf Example.} Consider the matrix and corresponding SR graph shown in Figure~\ref{exallo}. Observation of the SR graph tells us that the determinant expansion of the matrix contains two terms, and these are of the same sign. Details are provided in this figure legend. 

\begin{figure}[h]
\begin{minipage}{0.4\textwidth}
\[
M = \left[
\begin{array}{rrr}
a & b & \,\,\,\,\,0\\-c &  0 & \,\,d\\ 0 &  -e & \,\,f
\end{array}
\right]
\]

\end{minipage}
\hfill
\begin{minipage}{0.55\textwidth}
\begin{tikzpicture}[domain=0:4,scale=0.65]

\node at (1,4) {$S_1$};
\node at (1,1) {$R_1$};
\node at (4,4) {$R_2$};
\node at (4,1) {$S_2$};
\node at (7,4) {$S_3$};
\node at (7,1) {$R_3$};

\draw[-, line width=0.04cm] (1,1.5) -- (1, 3.5);
\draw[-, line width=0.04cm] (7,1.5) -- (7, 3.5);

\draw[-, dashed, line width=0.04cm] (1.5,1) -- (3.5, 1);
\draw[-, line width=0.04cm] (4.5,1) -- (6.5, 1);
\draw[-, line width=0.04cm] (1.5,4) -- (3.5, 4);
\draw[-, dashed, line width=0.04cm] (4.5,4) -- (6.5, 4);

\node at (2.5, 4.3) {$\scriptstyle{b}$};
\node at (2.5, 1.3) {$\scriptstyle{c}$};

\node at (5.5, 4.3) {$\scriptstyle{e}$};
\node at (5.5, 1.3) {$\scriptstyle{d}$};

\node at (1.3, 2.5) {$\scriptstyle{a}$};
\node at (7.3, 2.5) {$\scriptstyle{f}$};

\end{tikzpicture}
\end{minipage}
\caption{\label{exallo} {\em Left.} A matrix $M$. Assume that $a,b,c,d,e,f > 0$. {\em Right.} The corresponding SR graph. Identifying edges in the graph with their labels, the graph contains two term subgraphs, $E_1 = \{a, d, e\}$ and  $E_2 = \{b, c, f\}$. As $E_1 \cup E_2$ contains a single o-cycle $[a, b, e, f, d, c]$, we can infer from Lemma~\ref{allocycles} that the corresponding terms $T_1$ and $T_2$ in the expansion of $\mathrm{det}(M)$ satisfy $\mathrm{sign}(T_1) = \mathrm{sign}(T_2)$. Indeed it is easy to confirm that $\mathrm{det}(M) = T_1 + T_2$ where $T_1 = ade$ and $T_2 = bcf$.}
\end{figure}

The next lemma shows that having cycles which are both e-cycles and s-cycles in a graph means that some terms in a determinant expansion sum to zero. 
\begin{lem}
\label{oneeands}
Consider any square submatrix $S(\gamma|\delta)$ of a matrix $S$. Let $\alpha$ and $\beta$ be permutations of $\delta$ such that $T_\alpha$ and $T_\beta$ are nonzero terms in the determinant expansion of $S[\gamma|\delta]$. Assume that $E_\alpha \cup E_\beta$ contains exactly one cycle $C$, and this cycle is both an e-cycle and an s-cycle. Then $T_\alpha + T_\beta = 0$. 
\end{lem}
\begin{pf}
Let $|\gamma| = |\delta| = k$. By definition
\[
T_\alpha + T_\beta = P(\alpha)\prod_{i =1}^k S_{\gamma_i \alpha_i} + P(\beta)\prod_{i =1}^k S_{\gamma_i \beta_i}
\]
As usual, let $\theta$ be the set of indices for which $\alpha_i \not = \beta_i$ so that by assumption, $S_{\gamma_i \alpha_i}, S_{\gamma_i \beta_i} \in C$ iff $i \in \theta$. Defining $C_1 = \{S_{\gamma_i \alpha_i}\}_{i \in \theta}$ and $C_2 = \{S_{\gamma_i \beta_i}\}_{i \in \theta}$ gives us a disconnecting partition of $C$. We can write
\begin{eqnarray}
\hspace{-1cm}T_\alpha + T_\beta &= &\left(\prod_{i \in \{1, \ldots, k\}\backslash\theta} S_{\gamma_i \alpha_i}\right)\left(P(\alpha)\prod_{i \in \theta} S_{\gamma_i \alpha_i} + P(\beta)\prod_{i \in \theta} S_{\gamma_i \beta_i}\right) \nonumber \\
& = & P(\alpha)\left(\prod_{i \in \{1, \ldots, k\}\backslash\theta} S_{\gamma_i \alpha_i}\right)\left(\prod_{i \in \theta} S_{\gamma_i \alpha_i} + P(\beta\circ\alpha^{-1})\prod_{i \in \theta} S_{\gamma_i \beta_i}\right) \nonumber \\
& = & P(\alpha)\left(\prod_{i \in \{1, \ldots, k\}\backslash\theta} S_{\gamma_i \alpha_i}\right)\left(\mathrm{sign}(C_1)\mathrm{val}(C_1) + P(\beta\circ\alpha^{-1})\mathrm{sign}(C_2)\mathrm{val}(C_2)\right)\,. \nonumber 
\end{eqnarray}
$\beta \circ \alpha^{-1}$ can be written as a single cycle of length $|\theta|$, and so from Lemma~\ref{permsigns}, $P(\beta \circ \alpha^{-1}) = (-1)^{|\theta|-1}$. I.e.,
\[
T_\alpha + T_\beta = P(\alpha)\left(\prod_{i \in \{1, \ldots, k\}\backslash\theta} S_{\gamma_i \alpha_i}\right)\left(\mathrm{sign}(C_1)\mathrm{val}(C_1) + (-1)^{|\theta|-1}\mathrm{sign}(C_2)\mathrm{val}(C_2)\right)
\]
Since $C$ is an e-cycle we have
\[
\mathrm{sign}(C_2)/\mathrm{sign}(C_1) = \mathrm{sign}(C_1)\mathrm{sign}(C_2) = \mathrm{sign}(C) = (-1)^{|\theta|}.
\] 
Substituting into the expression for $T_\alpha + T_\beta$, we get:
\[
T_\alpha + T_\beta = P(\alpha)\left(\prod_{i \in \{1, \ldots, k\}\backslash\theta} S_{\gamma_i \alpha_i}\right)\mathrm{sign}(C_1)\left(\mathrm{val}(C_1) - \mathrm{val}(C_2)\right)\,.
\]
However since $C$ is an s-cycle, $\mathrm{val}(C_1) - \mathrm{val}(C_2)=0$, giving $T_\alpha + T_\beta = 0$.
\qquad \qed
\end{pf}

{\bf Example.} Consider the matrix and corresponding SR graph shown in Figure~\ref{exes}. Observation of the SR graph tells us that the matrix is singular. Details are in the figure legend.

\begin{figure}[h]
\begin{minipage}{0.4\textwidth}
\[
M = \left[
\begin{array}{rrr}
-a & b & \,\,\,\,\,0\\-c &  0 & \,\,b\\ 0 &  -c & \,\,a
\end{array}
\right]
\]

\end{minipage}
\hfill
\begin{minipage}{0.55\textwidth}
\begin{tikzpicture}[domain=0:4,scale=0.65]

\node at (1,4) {$S_1$};
\node at (1,1) {$R_1$};
\node at (4,4) {$R_2$};
\node at (4,1) {$S_2$};
\node at (7,4) {$S_3$};
\node at (7,1) {$R_3$};

\draw[-, dashed, line width=0.04cm] (1,1.5) -- (1, 3.5);
\draw[-, line width=0.04cm] (7,1.5) -- (7, 3.5);

\draw[-, dashed, line width=0.04cm] (1.5,1) -- (3.5, 1);
\draw[-, line width=0.04cm] (4.5,1) -- (6.5, 1);
\draw[-, line width=0.04cm] (1.5,4) -- (3.5, 4);
\draw[-, dashed, line width=0.04cm] (4.5,4) -- (6.5, 4);

\node at (2.5, 4.3) {$\scriptstyle{b}$};
\node at (2.5, 1.3) {$\scriptstyle{c}$};

\node at (5.5, 4.3) {$\scriptstyle{c}$};
\node at (5.5, 1.3) {$\scriptstyle{b}$};

\node at (1.3, 2.5) {$\scriptstyle{a}$};
\node at (7.3, 2.5) {$\scriptstyle{a}$};

\end{tikzpicture}
\end{minipage}
\caption{\label{exes} {\em Left.} A matrix $M$. Assume that $a, b, c > 0$. {\em Right.} The corresponding SR graph. There are two term subgraphs $E_1 = \{S_1\!-\!R_2,\, S_3\!-\!R_3,\, S_2\!-\!R_1\}$ and  $E_2 = \{S_1\!-\!R_1,\, S_3\!-\!R_2,\, S_2\!-\!R_3\}$. As $E_1 \cup E_2$ contains a single e-cycle which can be computed to be an s-cycle, Lemma~\ref{oneeands} implies that that the corresponding terms $T_1$ and $T_2$ in the expansion of $\mathrm{det}(M)$ satisfy $T_1 + T_2 = 0$, and hence that $M$ is singular. Indeed it is easy to confirm that $\mathrm{det}(M) = T_1 + T_2$ where $T_1 = abc$ and $T_2 = -abc$.}
\end{figure}

The next lemma is a consequence of the fact that if a matrix is nonsingular, then not all terms can be paired off in the way carried out in Lemma~\ref{oneeands}. 
\begin{lem}
\label{nonsing}
Consider any nonsingular square submatrix $S(\gamma|\delta)$ of a matrix $S$. Let $\alpha$ and $\beta$ be permutations of $\delta$ such that $T_\alpha$ and $T_\beta$ are nonzero terms in the determinant expansion of $S[\gamma|\delta]$. Assume that $E_\alpha \cup E_\beta$ contains exactly one cycle $C$, which is both an e-cycle and an s-cycle. Define $C_1 = C \cap E_\alpha$ and $C_2 = C \cap E_\beta$ so that $\{C_1, C_2\}$ is a disconnecting partition of $C$. Then $S[\gamma|\delta]$ must contain a term $T_\sigma$ such that $C_1 \not \subset E_\sigma$ and $C_2 \not \subset E_\sigma$.
\end{lem}
\begin{pf}
Assume the contrary. Take any term $T_\sigma$ in $S[\gamma|\delta]$. If $E_\sigma$ contains all the edges from $C_1$, then we can construct a new term subgraph $E_{\tau} = (E_\sigma\backslash C_1) \cup C_2$ and a corresponding term $T_\tau$ in $S[\gamma|\delta]$. Alternatively if $E_\sigma$ contains all the edges from $C_2$, define $E_{\tau} = (E_\sigma\backslash C_2) \cup C_1$ with corresponding term $T_\tau$. By construction, $E_\sigma \cup E_\tau$ contains only one cycle which is an e-cycle and an s-cycle and so, by Lemma~\ref{oneeands}, $T_\sigma+T_\tau = 0$. Thus all terms in $S[\gamma|\delta]$ fall into pairs which sum to zero and $S(\gamma|\delta)$ is singular.
\qquad \qed
\end{pf}

The next lemma tells us a fact which is geometrically obvious about how a term subgraph can intersect a cycle: either it contains all members in one half of a disconnecting partition of the cycle, or it contains an edge not in the cycle but incident on a vertex in the cycle. 
\begin{lem}
\label{termsets}
Consider a square matrix $S$ and the associated SR graph $G$. Assume that $G$ contains a cycle $C$. Let $\{C_1, C_2\}$ be a disconnecting partition of $C$, and $E$ be any term subgraph in $G$. If $C_1 \not \subset E$ and $C_2 \not \subset E$, then $E$ contains an edge incident on an S-vertex in $C$ but not itself in $C$.
\end{lem}
\begin{pf}
Let $|C| = 2k$. Since $E$ is a term subgraph in $G$, it has an edge incident on each  vertex in $C$. Either some of these are not in $C$, or $E$ contains $k$ edges in $C$. In the latter case, either they are all in $C_1$, or they are all in $C_2$, or two of them share a vertex. The first two possibilities are ruled out by assumption, and since $E$ is a term subgraph, the third situation is not possible. So $E$ must contain an edge $e$ incident on a vertex in $C$, but not itself lying in $C$. If this is an S-vertex we are done. If it is an R-vertex, but $e$ is also incident on another vertex in $C$, then this second vertex must be an S-vertex, and again we are done. Otherwise $e$ is incident on an R-vertex in $C$ and some vertex outside $C$, in which case $E\backslash\{e\}$ contains edges incident on $k$ S-vertices and $k-1$ R-vertices in $C$. There must hence be some edge $e' \in E\backslash\{e\}$ incident on an S-vertex in $C$ but not itself in $C$.
\qquad \qed
\end{pf}

Note that in the statement and proof of the above lemma we could exchange ``S-vertex'' and ``R-vertex''.

Our final preliminary lemma tells us that if an S-to-R path ``slices'' an e-cycle in a particular way, this implies the existence of two e-cycles with S-to-R intersection. 
\begin{lem}
\label{ecycles}
Consider a matrix $S$ and the associated SR graph $G$. Assume that $G$ contains an e-cycle $C$. Assume that there is a path $D$ in $G$ joining an S-vertex in $C$ to an R-vertex in $C$, but such that $D$ and $C$ are edge-disjoint. Then $G$ contains two e-cycles with S-to-R intersection. 
\end{lem}

\begin{pf}
Let $C$ be of length $2r$ with $C = \{S_{a_1 b_1}, S_{a_1 b_2}, S_{a_2 b_2}, \ldots, S_{a_r b_r}, S_{a_r b_1}\}$, and let $D$ join S-vertex $a_j$ to R-vertex $b_k$. Decompose $C$ into the two edge-disjoint paths $C = C_1 \cup C_2$ where
\[
C_1 = \{S_{a_j b_j}, S_{a_{j+1} b_j}, \ldots, S_{a_k b_k}\},  \quad C_2 = \{S_{a_{k+1} b_k}, S_{a_{k+1} b_{k+1}}, \ldots, S_{a_j b_{j-1}}\}\,.
\]
Note that $|C_1|$ and $|C_2|$ are both odd since $C_1$ and $C_2$ are S-to-R paths. The situation is illustrated in Figure~\ref{paths}. 

\begin{figure}[h]
\begin{minipage}[h]{\textwidth}
\begin{center}
\includegraphics[width=0.4\textwidth]{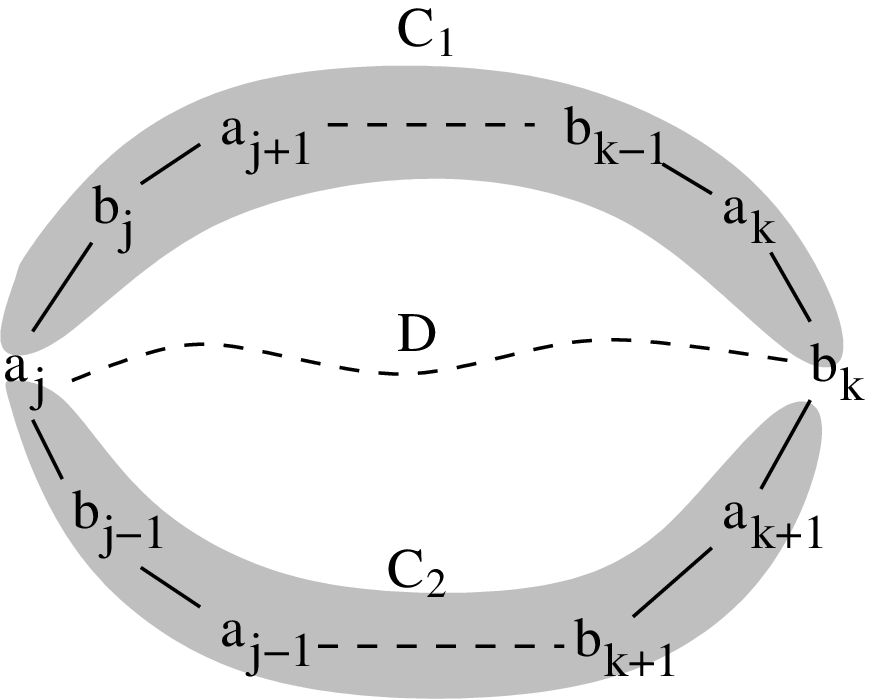}
\end{center}
\end{minipage}
\begin{minipage}[h]{\textwidth}
\caption{\label{paths}A schematic representation of the situation in Figure~\ref{ecycles}. $C_1$ and $C_2$ are S-to-R paths which partition $C$. $D$ joins an S-vertex to an R-vertex in $C$ but is edge-disjoint from $C$. (Note that in this figure, bold lines indicate edges of arbitrary sign, while dashed lines indicate paths of arbitrary length.)}
\end{minipage}
\end{figure}

There are two cycles $D_1 \equiv D \cup C_1$ and $D_2 \equiv D \cup C_2$ such that $D_1 \cap D_2 = D$, $D_1 \cap C = C_1$ and $D_2 \cap C = C_2$. We will show that $D_1$ and $D_2$ have opposite parity and hence one of them must be an e-cycle. Note that $|C_1| + |C_2| = 2r$, $|D_1| = |D| + |C_1|$ and $|D_2| = |D| + |C_2|$. Subtracting the last two expressions gives $|D_2| -|D_1| = |C_2| - |C_1| = 2(r - |C_1|)$. 

The fact that $C$ is an e-cycle means that either 
\begin{enumerate}
\item $r$ is odd and $\mathrm{sign}(C) = -1$, or
\item $r$ is even and $\mathrm{sign}(C)= +1$. 
\end{enumerate}

{\bf Case 1.} Since $\mathrm{sign}(C) = \mathrm{sign}(C_1)\,\mathrm{sign}(C_2) = -1$:
\[
\mathrm{sign}(C_1) = -\mathrm{sign}(C_2)
\]
Now $r$ and $|C_1|$ are odd, so $r - |C_1|$ is even, so $|D_2| - |D_1|= 2(r - |C_1|)$ is a multiple of $4$, i.e., $(-1)^{|D_2|/2} = (-1)^{|D_1|/2} $. Then 
\[
P(D_1) = (-1)^{|D_1|/2}\mathrm{sign}(C_1)\mathrm{sign}(D) = (-1)^{|D_2|/2}[-\mathrm{sign}(C_2)]\mathrm{sign}(D) = -P(D_2).
\]
So $P(D_2) = -P(D_1)$ and one of $D_1$ or $D_2$ must be an e-cycle. 

{\bf Case 2.} This time $\mathrm{sign}(C) = \mathrm{sign}(C_1)\,\mathrm{sign}(C_2) = 1$ so
\[
\mathrm{sign}(C_1) = \mathrm{sign}(C_2)
\]
Now $r$ is even and $|C_1|$ is odd, so $r - |C_1|$ is odd. As a result $|D_2| - |D_1| = 2(r - |C_1|)$ is not a multiple of $4$, and $(-1)^{|D_2|/2} = (-1)^{|D_1|/2 + 1} $. Again we get $P(D_2) = -P(D_1)$ and one of $D_1$ or $D_2$ must be an e-cycle. 

In each case, one of $D_1$ or $D_2$ is an e-cycle. Moreover both $D_1$ and $D_2$ intersect $C$ along an S-to-R path (either $C_1$ or $C_2$). As this is the unique component of their intersection, they have S-to-R intersection. 
\qquad \qed
\end{pf}

\section{Relationship between sign nonsingularity and o-cycles}

There is a very simple and elegant relationship between properties of square submatrices of a matrix and the non-existence of e-cycles in the corresponding SR graph. The results we prove in this section are weaker than our main result in the next section, but have a certain generality to them. At an abstract level, the results simply describe an interesting relationship between a matrix and an associated bipartite graph. Comments in Chapter~3 of \cite{brualdi} suggest that these results may be known or suspected, but we have not found a proof in the literature. In any case they are an easy corollary of the previous lemmas. First, analogous to the definition of sign-nonsingularity, we define a square matrix $S$ to be {\bf sign singular}, if all matrices with the same sign pattern as $S$ are singular. 

\begin{thm}
\label{SNS}
Consider a matrix $S$ and the associated SR graph $G$. If all cycles in $G$ are o-cycles, then all square submatrices of $S$ are either sign nonsingular or sign singular. 
\end{thm}

\begin{pf}
Consider any square submatrix $S(\gamma|\delta)$ of $S$ with rows and columns indexed in the usual way by ordered sets $\gamma$ and $\delta$ respectively, where $|\gamma| = |\delta| = k$. Consider permutations of $\delta$ and the corresponding terms in $S[\gamma|\delta]$:

1) If all such terms are zero then the $S(\gamma|\delta)$ is sign singular. 

2) If there is a single $\alpha$ for which $T_\alpha$ is nonzero then $S[\gamma|\delta]= T_\alpha$ so clearly $S(\gamma|\delta)$ is sign nonsingular.

3) Consider any pair of nonzero terms in $S[\gamma|\delta]$ corresponding to permutations $\alpha$ and $\beta$ of $\delta$. Since all cycles are odd, Lemma~\ref{allocycles} gives us that $\mathrm{sign}(T_\alpha) = \mathrm{sign}(T_\beta)$. Since $\alpha$ and $\beta$ were arbitrary, this means that all nonzero terms in $S[\gamma|\delta]$ have the same sign. 
\qquad \qquad \qed
\end{pf}

The previous theorem has a converse:

\begin{thm}
\label{SNS1}
Consider a matrix $S$ and the associated SR graph $G$. If all square submatrices of $S$ are either sign nonsingular or sign singular, then all cycles in $G$ are o-cycles. 
\end{thm}

\begin{pf}
Assume that $G$ has an $e$-cycle $C$ of length $2r$ including S-vertices from a set $\gamma=\{\gamma_1, \gamma_2, \ldots, \gamma_r\}$ and R-vertices from a set $\delta = \{\delta_1, \delta_2, \ldots, \delta_r\}$. There is some permutation $\alpha$ of $\delta$ such that $C$ consists of the edges 
\[
\bigcup_{i\in \{1, \ldots, r\}}\{S_{\gamma_i\alpha_i}, S_{\gamma_i\alpha_{i+1}}\}.
\]
Setting $\beta_i = \alpha_{i+1}$ defines a permutation $\beta$ of $\delta$. Clearly $S(\gamma|\delta)$ is not sign singular since
\[
T_\alpha = P(\alpha)\prod_{i = 1}^r S_{\gamma_i\alpha_i} \quad \mbox{and}\quad T_\beta = P(\beta)\prod_{i=1}^r S_{\gamma_i\beta_i}
\]
are nonzero terms in $S[\gamma|\delta]$. Since $E_\alpha \cup E_\beta = C$, and $C$ is an e-cycle, applying Eq.~(\ref{prodform}) gives:
\[
\mathrm{sign}(T_\alpha T_\beta)  = (-1)^{1} = -1\,.
\]
As $T_\alpha$ and $T_\beta$ have opposite signs, $S(\gamma|\delta)$ fails to be sign nonsingular. 
\qquad \qed
\end{pf}

Note that in the terminology of \cite{brualdi} a matrix which is either sign nonsingular or sign singular is a matrix with {\bf signed determinant}. We can state the previous two theorems as a single result:
\begin{cor}
\label{SNS2}
Consider a matrix $S$ and the associated SR graph $G$. The following two statements are equivalent:
\begin{enumerate}
\item All square submatrices of $S$ have signed determinant.
\item All cycles in $G$ are o-cycles. 
\end{enumerate}
\end{cor}
\begin{pf}
This follows immediately from Theorems~\ref{SNS}~and~\ref{SNS1}.
\qquad \qed
\end{pf}

{\bf Remark.} Corollary~\ref{SNS2} is a general statement about matrices and the corresponding bipartite graphs. Note that insisting that all square submatrices of a matrix have signed determinant is considerably more restrictive than insisting that they are all either sign nonsingular or singular (i.e. that the matrix is SSD). This is because it is a frequent occurrence for a square matrix to be singular while failing to be {\em sign singular}.

\section{A graph-theoretic condition ensuring injectivity}

Define the following condition on any SR graph:

{\bf Condition ($*$)}: all e-cycles in the SR graph are s-cycles, and no two e-cycles have S-to-R intersection. 

In \cite{craciun1} it was shown that Condition~($*$) on the SR graph of a reaction system with mass-action kinetics and with outflows is sufficient to ensure injectivity of the system. For {\bf N1C} reaction systems this means, by results in \cite{banajiSIAM}, that Condition~($*$) is sufficient to ensure that the stoichiometric matrix $S$ of the system is ``weakly sign determined'' (WSD), a less restrictive condition on $S$ than the requirement that it should be SSD. However results in \cite{craciun1} and \cite{banajiSIAM} left open the question of whether Condition~($*$) implies that the stoichiometric matrix is SSD giving injectivity for essentially arbitrary kinetics. Below we show that this is indeed the case -- Condition ($*$) implies that the stoichiometric matrix is SSD. Since SSD implies WSD, as a corollary we reproduce the result that can be inferred from \cite{craciun1} and \cite{banajiSIAM}. We also show by example that Condition~($*$) is not a necessary condition for a matrix to be SSD.

\subsection{Condition ($*$) ensures that a matrix is SSD}
\label{secsuff}

Theorem~\ref{suff} is the main result of this paper: that an {\bf N1C} reaction system whose SR graph fulfils Condition~($*$) has stoichiometric matrix which is SSD, and is hence, with the outflow conditions detailed at the beginning, incapable of multiple equilibria.

{\bf PROOF of Theorem~\ref{suff}.} Assume that $S$ is not SSD, i.e., there exists some square submatrix $S(\gamma|\delta)$ of $S$ which is neither sign nonsingular nor singular. Since $1 \times 1$ matrices are trivially sign nonsingular or singular, $|\gamma| = |\delta| \geq 2$. The proof will proceed by showing that the corresponding subgraph $G(\gamma|\delta)$ either contains an e-cycle which fails to be an s-cycle or contains two e-cycles which have S-to-R intersection. 

Firstly, it is immediate from Theorem~\ref{SNS} that $G(\gamma|\delta)$ contains an e-cycle. If $|\gamma| = 2$, then there are exactly two terms in $S[\gamma|\delta]$, i.e., $S_{\gamma_1\delta_1}S_{\gamma_2\delta_2}$ and  $-S_{\gamma_1\delta_2}S_{\gamma_2\delta_1}$. Since $S(\gamma|\delta)$ is not sign nonsingular, 
\[\mathrm{sign}(S_{\gamma_1\delta_1}S_{\gamma_2\delta_2}) = \mathrm{sign}(S_{\gamma_1\delta_2}S_{\gamma_2\delta_1})\,.\]
Since $S(\gamma|\delta)$ is not singular, $S[\gamma|\delta] = S_{\gamma_1\delta_1}S_{\gamma_2\delta_2}-S_{\gamma_1\delta_2}S_{\gamma_2\delta_1} \not = 0$. Hence, by the definition of an s-cycle, the e-cycle in $G(\gamma|\delta)$ consisting of the edges $\{S_{\gamma_1\delta_1}, S_{\gamma_1\delta_2}, S_{\gamma_2\delta_2}, S_{\gamma_2\delta_1}\}$ is not an s-cycle and Condition ($*$) is violated.

So now assume $|\gamma| \geq 3$. Consider two terms $T_\alpha$ and $T_{\eta}$ of opposite sign in $S[\gamma|\delta]$ corresponding to subgraphs $E_\alpha$ and $E_{\eta}$ in $G$. By Lemma~\ref{allocycles}, $E_\alpha \cup E_{\eta}$ contains an e-cycle. Take an e-cycle $C\subset (E_\alpha \cup E_{\eta})$ and construct a disconnecting partition of $C$ into $C_1 = C \cap E_\alpha$ and $C_2 = C \cap E_\eta$ according to Eq.~(\ref{cycdecomp}). Let $|C| = 2k$ for some $2 \leq k \leq |\gamma|$ so that $|C_1| = |C_2| = k$. Define $E_\beta = (E_\alpha\backslash C_1) \cup C_2$. $E_\beta$ has been constructed so $E_\alpha \cup E_\beta$ contains only one cycle (i.e., $C$). Corresponding to $E_\alpha$ and $E_\beta$ are nonzero terms $T_\alpha$ and $T_\beta$ in $S[\gamma|\delta]$. Bearing in mind that there is only one cycle in $E_\alpha \cup E_\beta$ and this is even, Eq.~(\ref{prodform}) gives: 
\[
\mathrm{sign}(T_\alpha T_\beta)  = (-1)^{1} = -1.
\]

There are now two cases to consider. If $T_\alpha + T_\beta \not = 0$, then $\mathrm{val}(E_\alpha) \not = \mathrm{val}(E_\beta)$. I.e., $\mathrm{val}((E_\alpha \backslash C) \cup (E_\alpha \cap C)) \not = \mathrm{val}((E_\beta \backslash C) \cup (E_\beta \cap C))$. Since $E_\alpha \backslash C = E_\beta \backslash C$, we get $\mathrm{val}(E_\alpha \cap C) \not = \mathrm{val}(E_\beta \cap C)$. Since $E_\alpha \cap C$ and $E_\beta \cap C$ together make up a disconnecting partition of $C$, this means that $C$ is not an s-cycle and we are done. 

So assume that  $T_\alpha + T_\beta = 0$. Now by Lemma~\ref{nonsing} we can find a permutation $\sigma$ of $\delta$, and corresponding term $T_\sigma$ in $S[\gamma|\delta]$ such that $C_1 \not \subset E_\sigma$ and $C_2 \not \subset E_\sigma$. As a result, by Lemma~\ref{termsets}, $E_\sigma$ contains some edge incident on an S-vertex in $C$, but not itself in $C$. In other words, there is an index $q$ such that edges $S_{\gamma_q\alpha_q}, S_{\gamma_q\beta_q} \in C$, but $\sigma_q \not = \alpha_q$ and $\sigma_q \not = \beta_q$ so that $S_{\gamma_q\sigma_q} \not \in C$.

Consider $\sigma\circ \alpha^{-1}$ as the product of disjoint cycles. Since $\sigma_q \not = \alpha_q$, a nontrivial cycle involves $\sigma_q$. Further, since $\sigma_q \not = \beta_q$, the corresponding cycle in $G$ is distinct from $C$ but intersects $C$ at S-vertex $\gamma_q$. Let this cycle be termed $C_{\alpha\sigma}$ (to remind us that it is composed of edges from $E_\alpha$ and $E_\sigma$) and have length $2r$ ($2 \leq r \leq |\gamma|$). In the usual way, we can follow $C_{\alpha\sigma}$: there are distinct indices $a_1, \ldots,  a_r \in \gamma$  and $b_1, \ldots, b_r \in \delta$ such that the edges $S_{a_j b_j}$ occur in $E_\alpha$ and $S_{a_j b_{j+1}}$ occur in $E_\sigma$. Assume (without loss of generality, i.e., by reordering the sets $a_i$ and $b_i$ if necessary) that $a_1 = \gamma_q$ and $b_2 = \sigma_q$, so that $S_{a_1 b_2}$ is the edge incident on an S-vertex in $C$, but not itself lying in $C$. 

Follow the cycle $C_{\alpha\sigma}$ starting at $S_{a_1 b_2} \in T_\sigma$ i.e.,
\[
S_{a_1 b_2} \in T_\sigma,\,\, S_{a_2 b_2} \in T_\alpha,\,\, S_{a_2 b_3} \in T_\sigma, \,\,\ldots
\]
Since this is a cycle, eventually some vertex from the sequence of alternating R- and S-vertices indexed by $(b_2, a_2, b_3, \ldots)$ must be a vertex in $C$. But this cannot first happen at an S-vertex. Suppose the contrary and S-vertex $a_j$ is in $C$, while R-vertex $b_j$ is not. Since edge $S_{a_j b_j}$ is in $T_\alpha$, but not in $C$, this implies that S-vertex $a_j$ has three edges from $T_\alpha \cup T_\beta$ incident on it  -- the two edges in $C$ along with the edge $S_{a_j b_j}$. But this is impossible from the discussion in Section~\ref{detperm}. This means that from the vertex sequence $(b_2, a_2, b_3, \ldots)$ the first vertex to lie in $C$ must be an R-vertex. Let this vertex be $b_j$  ($j=2$ is possible). Define the path $D \equiv \{S_{a_1 b_2}, S_{a_2 b_2}, S_{a_2 b_3}, \ldots, S_{a_{j-1}b_j}\}$. $D$ is an S-to-R path starting and terminating at vertices in $C$ but edge-disjoint from $C$. 

Now applying Lemma~\ref{ecycles}  shows that there are two e-cycles in $G$ (one of which is $C$ and one of which is made up of $D$ and some part of $C$) which have S-to-R intersection. Thus if the stoichiometric matrix $S$ is not SSD, then the associated graph $G$ necessarily fails Condition~($*$). The result is proved.
\qquad \qed

\subsection{Condition ($*$) is not necessary for SSD (or WSD)}
\label{secnecc}

We present an example to illustrate that Condition~($*$) on the SR graph is not necessary to give a system with an SSD matrix. Consider the system of three reactions

\begin{equation}
\label{StoRcounter}
D \rightleftharpoons A + B +C,\quad E \rightleftharpoons A + B + C, \quad F \rightleftharpoons A + B. 
\end{equation}

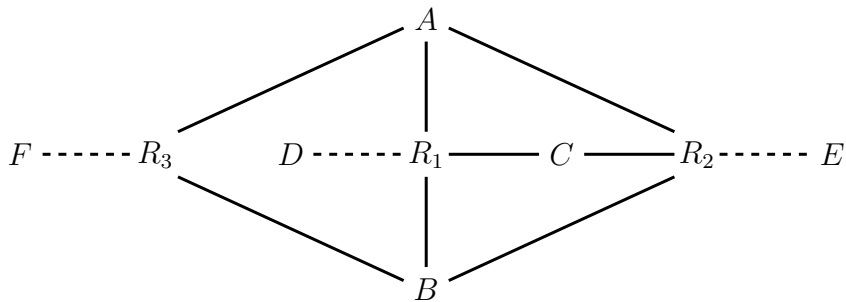
\begin{figure}[h]
\begin{minipage}{\textwidth}
\begin{center}
\begin{tikzpicture}[domain=0:4,scale=0.6]

\node at (1,4) {$F$};
\node at (4,4) {$R_3$};
\node at (7,4) {$D$};
\node at (10,4) {$R_1$};
\node at (13,4) {$C$};
\node at (16,4) {$R_2$};
\node at (19,4) {$E$};

\node at (10,1) {$B$};
\node at (10,7) {$A$};

\draw[-, dashed, line width=0.04cm] (1.5, 4) -- (3.5,4);
\draw[-, dashed, line width=0.04cm] (7.5, 4) -- (9.5,4);
\draw[-, line width=0.04cm] (10.5, 4) -- (12.5,4);
\draw[-, line width=0.04cm] (13.5, 4) -- (15.5,4);
\draw[-, dashed, line width=0.04cm] (16.5, 4) -- (18.5,4);

\draw[-, line width=0.04cm] (10, 1.5) -- (10,3.5);
\draw[-, line width=0.04cm] (10, 4.5) -- (10,6.5);

\draw[-, line width=0.04cm] (9.5, 1.2) -- (4.5,3.5);
\draw[-, line width=0.04cm] (10.5, 1.2) -- (15.5,3.5);

\draw[-, line width=0.04cm] (9.5, 6.8) -- (4.5,4.5);
\draw[-, line width=0.04cm] (10.5, 6.8) -- (15.5,4.5);

\end{tikzpicture}
\end{center}

\end{minipage}
\caption{\label{counter}The SR graph of reaction system~(\ref{StoRcounter}). Edge labels are all $1$ and have been omitted. As a consequence all cycles are s-cycles. However, there are a number of e-cycles with S-to-R intersection, for example the cycles $A\!-\!R_2\!-\!B\!-\!R_3$ and $A\!-\!R_2\!-\!C\!-\!R_1$ intersect along the S-to-R path $A\!-\!R_2$.}
\end{figure}

This system has stoichiometric matrix
\[
S = 
\left[
\begin{array}{rrr}
1 & 1 & 1\\1 & 1 & 1\\1 & 1 & 0\\-1 & 0 & 0 \\0 & -1 & 0\\0 & 0 & -1
\end{array}
\right]
\]
which can quickly be computed to be SSD, implying that the system with inflows and outflows forbids multiple equilibria. The SR graph for the system is shown in Figure~\ref{counter}. Although all e-cycles are s-cycles, there are e-cycles with S-to-R intersection, and thus the graph fails Condition~($*$). 

Thus Condition~($*$) is not necessary for injectivity in the case of a general system of {\bf N1C} reactions. Since the condition that $S$ is SSD is stronger than the condition that $S$ is WSD, clearly Condition~($*$) is not necessary for injectivity in the case of a system of mass action reactions. 

\section{Conclusions}
We have described several new results for the class of {\bf N1C} reaction systems with arbitrary chemical kinetics. The key has been to associate with any real matrix $S$, a signed, labelled, bipartite graph termed an SR graph. Given any matrix $S$ and associated graph $G$ we have shown that:
\begin{enumerate}
\item All cycles in $G$ are o-cycles iff all square submatrices of $S$ have signed determinant (Corollary~\ref{SNS2}).
\item Condition~($*$) on $G$ is sufficient, but not necessary, to guarantee that $S$ is SSD (Theorem~\ref{suff}, and the counterexample in Section~\ref{secnecc}).
\end{enumerate}
If $S$ is the stoichiometric matrix of a chemical reaction network, and we assume the outflow conditions in System~(\ref{basicout}), either graph-theoretic condition implies immediately that multiple equilibria cannot exist \cite{banajiSIAM}. Otherwise, they imply that multiple positive nondegenerate equilibria cannot exist \cite{craciun2}. These results apply to large classes of realistic biochemical reaction networks; for example, they apply to the enzymatic reaction networks discussed in \cite{craciun3}, without the assumption of mass-action kinetics. There are natural further extensions of this work to the situation where the {\bf N1C} condition is dispensed with. These topics are pursued in \cite{banajicraciun2}.

The identification of a condition on $G$ which is {\em equivalent} to $S$ being SSD remains an open -- and probably very difficult -- problem. However, there are certain ways in which $G$ can fail Condition~($*$) which ensure that $S$ fails to be SSD. These will be discussed in future work.

Finally, it should be mentioned that there are interesting relationships between the results in this paper and approaches to injectivity involving so-called ``interaction graphs''. The approach taken in \cite{soule} is most similar to that taken here, with the key difference that interaction graphs and SR graphs are formally quite different objects. The details of how interaction graph approaches and SR graph approaches are related to each other will be explored in future work. 

\section*{Acknowledgements}

MB's research was funded by EPSRC grant EP/D060982/1. GC acknowledges support from NIH grant R01GM086881. 

\bibliographystyle{elsart-num-sort}

\end{document}